\documentclass{amsart}

\newtheorem{theorem}{Theorem}[section]
\newtheorem{lemma}[theorem]{Lemma}

\theoremstyle{definition}

\newtheorem{proposition}[theorem]{Proposition}
\theoremstyle{remark}
\newtheorem{remark}[theorem]{Remark}

\numberwithin{equation}{section}

\begin{document}

\title{$\pi_1$ of Hamiltonian $S^1$ manifolds}

\author{Hui Li}
\address{Department of Mathematics, University of Illinois, 
Urbana-Champaign, IL, 61801} 
\curraddr{Department of Mathematics,
University of Illinois, Urbana-Champaign, IL, 61801}
\email{hli@math.uiuc.edu}


\date{January 10, 2002 and, in revised form, May 16, 2002.}


\begin{abstract}
   Let $(M, \omega)$ be a connected, compact symplectic manifold equipped with a Hamiltonian $S^1$ action.
  We prove that, as fundamental groups of topological spaces, $\pi_1(M)=\pi_1(\hbox{minimum})=\pi_1(\hbox{maximum})=\pi_1(M_{red})$,
 where $M_{red}$ is the symplectic quotient at any  value in the image of the moment map $\phi$.                                                                   
\end{abstract}

\maketitle

 Let $(M,\omega)$ be a connected, compact symplectic manifold equipped with a circle action. If the action is Hamiltonian, then the moment map
 $\phi: M\rightarrow \mathbb{R}$ is a perfect Bott-Morse function. Its critical sets are precisely the fixed point sets $M^{S^1}$ of the $S^1$ action, and
 $M^{S^1}$ is a disjoint union of symplectic submanifolds. Each fixed point set has even index. By \cite{A}, $\phi$ has a unique local minimum and a
 unique local maximum. We will use Morse theory to prove
 \begin{theorem}
 Let $(M^{2n}, \omega)$ be a connected, compact symplectic manifold equipped with a Hamiltonian $S^1$ action. Then, as fundamental groups
 of topological spaces, $\pi_1(M)=\pi_1(\hbox{minimum})=\pi_1(\hbox{maximum})=\pi_1(M_{red})$,
  where $M_{red}$ is the symplectic quotient at any  value in the image of the moment map $\phi$.
  \end{theorem}
  
  \begin{remark}\label{rem1}
    The theorem is not true for orbifold $\pi_1$ of $M_{red}$, as shown in the example below. (See \cite{C} or \cite{TY} for the definition of orbifold $\pi_1$).

   Let $a\in\hbox{im}(\phi)$, and $\phi^{-1}(a)=\{x\in M\mid \phi(x)=a\}$ be the level set. Define $M_a=\phi^{-1}(a)/S^1$ to be the symplectic quotient.

   Note that if $a$ is a regular value of $\phi$, and if the circle action on $\phi^{-1}(a)$ is not free,
   then $M_a$ is an orbifold, and we have an orbi-bundle:
   \begin{equation}\label{eq1} 
  \begin{array}{ccl}
   S^1 &\hookrightarrow & \phi^{-1}(a)\\
      &                & \downarrow\\
      &                 & M_a     
   \end{array}
\end{equation}
     
   If $a$ is a critical value of $\phi$, then $M_a$ is a stratified space.
   (\cite{SL}).

   Now, let $S^1$ act on $(S^2\times S^2, 2\rho\oplus\rho)$ (where $\rho$ is the standard symplectic
   form on $S^2$) by $\lambda (z_1, z_2)=(\lambda^2 z_1, \lambda z_2)$. Let $0$  be the minimal value of the moment map. Then for $a\in (1, 2)$, $M_a$
   is an orbifold which is homeomorphic to $S^2$ and has two $\mathbb{Z}_2$ singularities. The orbifold $\pi_1$ of $M_a$ is $\mathbb{Z}_2$, but the $\pi_1$ of $M_a$ as a 
   topological space is trivial.
  \end{remark}

  Let $a$ be a regular or a critical value of $\phi$. Define
  $$M^a=\{x\in M\mid \phi(x)\leq a\}.$$ 
  
  By Morse theory, we have the following lemmas about how $M^a$ and $\phi^{-1}(a)$ change when $\phi$ doesn't cross or crosses a critical level.
  \begin{lemma}\label{lem1}
  (Theorem 3.1 in \cite{M}) Assume $[a,b]\subset\hbox{im}(\phi)$ is an interval consisting of regular values, then $\phi^{-1}(a)$ is diffeomorphic to
   $\phi^{-1}(b)$.
   \end{lemma}
   
  \begin{lemma}\label{lem2} (See \cite{M} and \cite{B}) 
   Let $c\in (a,b)$ be the only critical value of $\phi$ in $[a,b]$, $F\subset \phi^{-1}(c)$ be the fixed point set component, $D^-$ be the negative
  disk bundle of $F$, and $S(D^-)$ be its sphere bundle. Then $M^b$ is homotopy equivalent to $M^a\cup_{S(D^-)}D^-$.
  \end{lemma}
  \begin{lemma}\label{lem3}
  Under the same hypothesis of Lemma~\ref{lem2}, $\phi^{-1}(a)\cup_{S(D^-)}D^-$ has the homotopy type of $\phi^{-1}(c)$.
  \end{lemma}
  
  \begin{proof}
   If $F$ is a point, then from  the  proof of Theorem 3.2 in \cite{M}, we can see that the region between $\phi^{-1}(a)\cup_{S(D^-)}D^-$ and $\phi^{-1}(c)$
   is homotopy equivalent to both $\phi^{-1}(a)\cup_{S(D^-)}D^-$ and $\phi^{-1}(c)$. (See $p_{18}$ and $p_{19}$ in \cite{M}).

    The same idea applies when $F$ is a submanifold. 
  \end{proof} 
  
   This lemma immediately implies the following        
  \begin{lemma}\label{lem4}
   Under the same hypothesis of Lemma~\ref{lem2}, $M_c$ has the homotopy type of $M_a\cup_{S(D^-)/S^1}D^-/S^1$.
  \end{lemma}

   We will also need
   \begin{lemma}\label{lem5}
    Assume $F$ is a critical set, $\phi(F)\in (a, b)$ and there are no other critical sets in $\phi^{-1}([a,b])$. If index(F)=2, then there is an 
   embedding $i$ from $F$ to $M_a$, such that $S(D^-)$ can be identified with the restriction of $\phi^{-1}(a)$ to $F$, i.e., we have the following bundle
   identification:
   \begin{equation}\label{eq2}
   \begin{array}{ccccc}
   S^1 & \hookrightarrow & S(D^-)   &\rightarrow &  \phi^{-1}(a)\\
       &                 & \downarrow &           &\downarrow \\
       &                 &F   & \overset{i}\longrightarrow &   M_a  \\
   \end{array}
   \end{equation}
   \end{lemma}
   \begin{proof}
     Assume that the positive normal bundle $D^+$ of $F$ has complex rank $m$. We may assume $\phi(F)=0$. By Lemma~\ref{lem1}, we can assume $a=-\epsilon$ and $b=+\epsilon$ 
    for $\epsilon$ small. By the equivariant symplectic embedding theorem (\cite{Ma}), a tubular neighborhood of $F$ is equivariantly diffeomorphic to 
    $P\times_{G}(\mathbb{C}\times\mathbb{C}^m)$,
   where $G=S^1\times U(m)$ and $P$ is a principal $G$-bundle over $F$.  
   And the moment map can be written $\phi=-p_0|z_0|^2+p_1|z_1|^2+...+p_m|z_m|^2$, where $p_0, p_1,..., p_m$ are positive integers. Then
   $\phi^{-1}(-\epsilon)=P\times_G(S^1\times\mathbb{C}^m)$, $M_{-\epsilon}=P\times_G(S^1\times\mathbb{C}^m)/S^1$. 
   $F=P\times_G(S^1\times 0)/S^1\subset M_{-\epsilon}$, and $S(D^-)=P\times_GS^1$ is the restriction of $\phi^{-1}(-\epsilon)$ to $F$.
 \end{proof}

     We are now ready to prove the theorem.
   \begin{proof}
   Let us put the critical values of $\phi$ in the order
   $$\hbox{minimal}=0<a_1<a_2<....<a_k=\hbox{maximal}.$$\\
   First, we prove $\pi_1(\hbox{minimum})=\pi_1(M_{red})$.

   For $a\in (0,a_1)$, by the equivariant symplectic embedding theorem, $\phi^{-1}(a)$ is a sphere bundle over the minimum. Assume the fiber 
   of this sphere bundle is $S^{2l+1}$ with $l\geq 0$. Then $M_a$ is diffeomorphic to a weighted $\mathbb{C}P^l$ bundle
   over the minimum (possibly an orbifold). The weighted $\mathbb{C}P^l$ is the symplectic reduction of $S^{2l+1}$ by the $S^1$ action with different weights.
   We can easily see that $S^{2l+1}\rightarrow $ weighted $\mathbb{C}P^l$ induces a surjection in $\pi_1$ since
   the inverse image of each point is connected. So the weighted $\mathbb{C}P^l$ is simply connected, hence $\pi_1(M_a)=\pi_1(\hbox{minimum})$.

   Next, let $b\in (a_1,a_2)$, and $F\subset \phi^{-1}(a_1)$ be the critical set. (If there are other critical sets on the same level, argue similarly
   for each connected component).

    By Lemma~\ref{lem4}, and the Van-Kampen theorem, we have
    $$\pi_1(M_{a_1})=\pi_1(M_a)*_{\pi_1(S(D^-)/S^1)}\pi_1(D^-/S^1)=\pi_1(M_a),$$
    since $S(D^-)/S^1$ is a weighted projectivized bundle over $F$, and $D^-/S^1$ is homotopy equivalent to $F$, so $\pi_1(S(D^-)/S^1)$ is isomorphic to
    $\pi_1(D^-/S^1)$.

    Similarly, using $-\phi$, we can obtain $\pi_1(M_b)=\pi_1(M_{a_1})$.

    By induction on the critical values,  and by repeating the argument each time $\phi$ crosses a critical level, we  see that if $a'\in (a_{k-1}, a_k)$,
   then $\pi_1(M_{a'})=\pi_1(\hbox{minimum})$. Similar to the proof of $\pi_1(M_a)=\pi_1(\hbox{minimum})$ when $a\in (0, a_1)$, we have 
    $\pi_1(M_{a'})=\pi_1(\hbox{maximum})$.

   Therefore we have proved that $\pi_1(M_{red})=\pi_1(\hbox{minimum})=\pi_1(\hbox{maximum})$.\\

   Next, we prove $\pi_1(M)=\pi_1(\hbox{minimum})$.

   Consider $M^a$, for $a\in (0, a_1)$. Since $M^a$ is a complex disk bundle over the minimum,  $\pi_1(M^a)=\pi_1(\hbox{minimum})=\pi_1(M_a)$.

    Consider  $b\in (a_1, a_2)$, and  let $F\subset \phi^{-1}(a_1)$ be the critical set.

     First assume index(F)=2. By Lemma~\ref{lem2}, and the Van-Kampen theorem, 
    $$\pi_1(M^b)=\pi_1(M^a)*_{\pi_1(S(D^-))}\pi_1(D^-)=\pi_1(M^a)*_{\pi_1(S(D^-))}\pi_1(F).$$
    Consider the fibration
   \begin{equation}\label{eq4}
   \begin{array}{ccl}
   S^1 &\hookrightarrow & S(D^-) \\
       &                & \downarrow \\
        &         &       F  \\
   \end{array}
   \end{equation}      
  and its homotopy exact sequence
  \begin{equation}\label{eq5}
  \cdot\cdot\cdot\rightarrow\pi_1(S^1)\overset{j}\rightarrow \pi_1(S(D^-))\overset{f}\rightarrow \pi_1(F)\rightarrow 0.
  \end{equation}

   The map $f$ is surjective. By Lemma~\ref{lem5}, the image of $\hbox{ker}(f)=\hbox{im}(j)$ in $\pi_1(M_a)$ is $0$. By induction, $\pi_1(M^a)=\pi_1(M_a)$. 
   So the image of $\hbox{ker}(f)$ in $\pi_1(M^a)$ is $0$.
   Hence $\pi_1(M^b)=\pi_1(M^a)=\pi_1(\hbox{minimum})$.

    If index(F)$>2$, then the corresponding map $\pi_1(S(D^-))\rightarrow\pi_1(F)$ is an isomorphism.  So we also have $\pi_1(M^b)=\pi_1(M^a)$.

    By induction, we see $\pi_1(M)=\pi_1(\hbox{minimum})$.
    \end{proof}
    
   \begin{remark}
   The proof that $\pi_1(\hbox{minimum})=\pi_1(M_{red})$ can be achieved by using known results about how the reduced space changes after $\phi$
   crosses a critical level. (See \cite{BP} for instance, or \cite{GS} where the action is semi-free). After the first induction step,
    when $\phi$ crosses a critical
   set $F$, if index(F)=2, then $M_a$ is homeomorphic to $M_{a_1}$; if index(F)$>2$, then $M_{a_1}$ can be obtained from $M_a$ by a blow-up
    followed by a blow-down. $M_b$ and $M_{a_1}$ are similarly related. Then we  modify the proof of D. McDuff's (\cite{Mc}) result
   \begin{proposition}
    If $\tilde{X}$ is the blow-up of $X$ along a submanifold $N$, then $\pi_1(\tilde{X})=\pi_1(X)$.
   \end{proposition}
  \end{remark}

   \subsubsection*{Acknowledgement}
    
   I would like to genuinely thank  Susan Tolman  for carefully reading this paper, for making corrections and for helping to improve the proof of 
   the theorem.

\bibliographystyle{amsplain}

\end{document}